\documentclass[10pt]{amsart}
\usepackage{amssymb, mathrsfs,amsmath}
\usepackage{latexsym}
\usepackage{amsfonts}
\usepackage{shadow}

\newtheorem{Pa}{Paper}[section]
\newtheorem{Tm}[Pa]{{\bf Theorem}}
\newtheorem{La}[Pa]{{\bf Lemma}}
\newtheorem{Cy}[Pa]{{\bf Corollary}}
\newtheorem{Ob}[Pa]{{\bf Observation}}

\newtheorem{Pn}[Pa]{{\bf Proposition}}

\newtheorem{Ex}[Pa]{{\bf Example}}

\def\e{\epsilon}
\def\C{\mathbb C}

\date{}
\author[D. Alpay]{Daniel Alpay}
\author[P. Jorgensen]{Palle Jorgensen}
\address{(DA) Department of Mathematics
\newline
Ben Gurion University of the Negev \newline P.O.B. 653,
\newline
Be'er Sheva 84105, \newline ISRAEL} \email{dany@math.bgu.ac.il}
\address{(PJ)
Department of Mathematics\newline 14 MLH \newline The University
of Iowa Iowa City,\newline IA 52242-1419 USA}
\email{jorgen@math.uiowa.edu}
\author[I. Lewkowicz]{Izchak Lewkowicz}
\address{(IL) Department of Electrical Engineering
\newline
Ben Gurion University of the Negev \newline P.O.B. 653,
\newline
Be'er Sheva 84105, \newline ISRAEL}
\email{izchak@ee.bgu.ac.il}
\thanks{This research is partially supported by the BSF
grant no. 2010117.}
\thanks{D. Alpay thanks the
Earl Katz family for endowing the chair which supported his
research.}
\title
[Parametrizations of all wavelet filters]
{
Parametrizations of all
wavelet filters: input-output and state-space}

\begin{document}
\begin{abstract}
We here use notions from the theory linear shift-invariant
dynamical systems to provide an explicit characterization,
both practical and computable, of all rational wavelet filters.
For a given $N$, $(N\geq 2)$ the number of inputs, the
construction is based on a factorization to an~ {\em elementary
wavelet filter} along with of $m$~ {\em elementary unitary}~
matrices. We shall call this $m$ the {\em index}~ of the filter.
It turns out that the resulting wavelet filter is of McMillan
degree $N\left(\frac{1}{2}(N-1)+m\right)$.

Moreover, beyond the parameters $N$ and $m$, one confine the
spectrum of the filters to lie in an open disk of radius $\rho$
(stable filters mean $\rho\in[0, 1]$ and for FIR take $\rho=0$).
Then all filters can be described by a convex set of parameters
$\left([0,~\pi)\times[0,~2\pi)^{2(N-1)}\times[0,~\rho)\right)^m$.

Rational wavelet filters bounded at infinity, admit state
space realization. The above input-output parametrization is
exploited for a step-by-step construction (where in each, the
index $m$ is increased by one) of state space model of
wavelet filters.
\end{abstract}

 \maketitle

\section{Introduction}
\setcounter{equation}{0}
\label{sec:1}
\subsection{Problem formulation - symmetries}
Over decades, filters have played a number of roles in both
signal processing and applied mathematics.
Of particular interest have been wavelet filters
\cite{BJ1}, \cite{BJ2}, \cite{J2}, \cite{HW},  \cite{Ma}.

To describe them let $\mathbb{T}$, $\mathbb{D}$ be the unit
circle, and the open unit disk, respectively, i.e.
\[
\mathbb{T}:=\{z\in\C~:~|z|=1~\}\quad\quad\quad
\mathbb{D}:=\{z\in\C~:~1>|z|~\}
\]
($\overline{\mathbb D}=\mathbb{D}\cup\mathbb{T}$ is the closed
unit disk). Roughly speaking wavelet filters are matrix valued
rational functions $Sy(z)$ and $An(z)$, analytic outside
$\mathbb{D}$, of a specific structure satisfying,
\begin{equation}\label{RE}
Sy(z)An(z)=I\quad\quad\quad\quad z\in\mathbb{T}.
\end{equation}
In signal processing terminology,
$An$ and $Sy$ are called the ``analysis" and ``synthesis" filters,
respectively, see e.g. \cite[Fig. 9.9]{SN}, \cite[Fig. 12.9-3]{Va},
\cite[Fig. 1]{XM}. The fact that $Sy$ is a (left) inverse of
$An$ is referred to as the ``perfect reconstruction" condition,
see e.g. \cite[Eq. (9.2)]{SN}, \cite[Section 5.6]{Va}.
\vskip 0.2cm

This work is aimed at three different communities: mathematicians
interested in classical analysis, signal processing engineers and
system and control engineers. Thus adopting the terminology
familiar to one audience, may intimidate or even alienate the
other.
For example what is known to engineers as McMillan degree also
arises in geometry of loop groups as an index. Filter banks and
multi-bands from signal processing turn into representations of
a Cuntz algebra, with the number of generating isometries in
the Cuntz algebra equal to the number of frequency bands.
%: ``the Cuntz algebras", ``multiresolution filter banks",
%and ``the system matrix associated with a state-space
%realization" are samples of terms, relevant to this work, from the
%jargon of (only) one of the above communities. 
Books like \cite{BN}, \cite{BJ2} and \cite{SN} and papers like
\cite{J1}, \cite{Wa} have made an impressive effort to be at
least ``bilingual". Lack of space prevents us from providing even a
concise dictionary of filter theory terms. Instead, we try to
employ only basic concepts or indicate for references
providing for the necessary background.
\vskip 0.2cm

Our paper is concerned with the use of time-frequency filters
in the construction of specific wavelets, for example in
$L^2(\mathbb R)$, or in in $L^2(\mathbb R^k)$. The use of these
filters offered a big boost to the list of constructive and
computable approaches. The significance of this lies both on
the theoretical side as well as the practical side; e.g., in
adaptive wavelet approximations; as well as the building of
fast computation of wavelet coefficients and reconstruction.
Signal and image processing have typically been the source of
time-frequency filters to be used in wavelets. Our present aim
is to supplement this with tools from systems theory; as well
as stressing the interdisciplinary unity of these themes. This
means we here focus on studying (wavelet) {\em filters},
leaving for a future work
the detailed wavelet analysis going into building
%wavelets
bases from filters. There are several reasons for our choice
of emphasis. Firstly, the process of building wavelets from
filters is covered in the literature; see e.g., \cite{BN, BJ0,
BJ1, BJ2}; but secondly, the analysis involved in this part of
the subject is a separate endeavor. Indeed, it offers a variety
of tools depending on the focus, for example an emphasis on the
theory, as opposed to ``tailoring" the basis to a specific
application.
%one of numerous applications.
We
have omitted a detailed account of it as it would take us too far
afield. The following sample references should help:
\cite{BFMP09, BlKr12, GlDA09,  Gu11,   Ha11,
LaRa06,  MiSu07, OlOl10, Rae09,SuVi10, YaZh08}.
\vskip 0.2cm

To simplify presentation we next {\em define}~ wavelet filters by
two symmetries: \eqref{eq:unit} and \eqref{eq:m1}.
\vskip 0.2cm

Denote by $\mathcal{U}_N$, with $N$ natural, the set of
$N\times N$-valued rational functions $U(z)$
unitary on the unit circle $\mathbb T$, i.e.
\begin{equation}\label{eq:unit}
\mathcal{U}_N=\{ U(z)~:~
U(z)^*U(z)=I_N\quad z\in\mathbb T~\}.
\end{equation}
Where $U^*$ denotes the complex conjugate transpose of $U$.
Whenever clear from the context we shall omit the subscript $N$
and simply write $\mathcal{U}$. The infinite-dimensional group
$\mathcal{U}$ in \eqref{eq:unit}, is well known in literature,
e.g. \cite{AlGo1} Versions of it are studied in physics under
the name loop-group \cite{DJ2}. In parts of the signal processing
literature $\mathcal{U}$ is referred to as the set of~
{\em paraunitary}~ matrices, see e.g. \cite[Section 5.1]{SN},
\cite[Chapter 6]{Va}. In this work we confine the discussion to
the subset of $\mathcal{U}$ which is analytic outside
$\mathbb{D}$, the open unit disk, (``Schur asymptotically
stable" in some communities). In electrical networks terminology,
it is sometimes referred to as ``lossless", see e.g. \cite[Sections
3.5 and 14.2]{Va}.
% With a slight abuse of notation, it will be denoted hereafter
% by $\mathcal{U}$. 
Generalizations of $\mathcal{U}$ to the unstable case, in
conjunction of wavelet filters, are addressed in
\cite{AJL}, \cite{AJLM}. See also Remark 4
in Section \ref{sec:10}.
\vskip 0.2cm

For $N\geq 2$, we shall find it convenient to denote by $\e$ the
following $N$-th root of unity,
\begin{equation}\label{eq:e}
\e:=e^{i\frac{2\pi}{N}}.
\end{equation}
We shall say that an $N\times N$-valued $(N\geq 2)$
rational function
$F(z)$ belongs to $\mathcal{C}_N$ if it is of the form,
\begin{equation}\label{eq:m1}
F(z)=\begin{pmatrix} \hat{f}_0(z)&
\hat{f}_0(\e z)&\cdots &\hat{f}_0(\e^{N-1}z)\\
\hat{f}_1(z)&\hat{f}_1(\e z)&\cdots &\hat{f}_1(\e^{N-1}z)\\
\vdots& &  & \\
\hat{f}_{N-1}(z)&\hat{f}_{N-1}(\e z)&\cdots &
\hat{f}_{N-1}(\e^{N-1}z)
\end{pmatrix},
\end{equation}
with $\e$ in \eqref{eq:e}. In signal processing literature
$F(z)$ in \eqref{eq:m1} is called a ``filter bank" e.g.
\cite[Chapter 4]{SN}, \cite[Chapter 11]{Va}. As before,
whenever clear from the context we shall
omit the subscript $N$ and simply write $\mathcal{C}$.
\vskip 0.2cm

In this work we focus our attention on wavelet filters, denoted
by $\mathcal{W}_N$, that is functions within $\mathcal{U}$
\eqref{eq:unit} which satisfy \eqref{eq:m1}, i.e.
\begin{equation}\label{eq:WaveletFilters}
\mathcal{W}_N:=\mathcal{U}_N\cap\mathcal{C}_N~.
\end{equation}
In terms of \eqref{RE} one simply takes,
\[
An(z)=W(z)\quad\quad {\rm and}\quad\quad
Sy(z)=W\left(\frac{1}{z^*}\right)^*,
\]
with $W\in\mathcal{W}$, see e.g. \cite[Fig. 4.3-12]{Va} (recall
that $z=\frac{1}{z^*}$ on $\mathbb{T}$). Up to Remarks 6, 7 in
Section \ref{sec:10},
we shall not relate to \eqref{RE} and only to
\eqref{eq:WaveletFilters}. The fact that wavelet filters satisfy
\eqref{eq:WaveletFilters} has been long recognized, see e.g.
\cite{Be}. In the spirit of \cite{BJ2}, we here use
\eqref{eq:WaveletFilters} as their characterization.
In signal processing literature, $N$ is referred to as ``the
number of bands of filter".
\vskip 0.2cm

The structure of this family suggests that its description is not
straightforward: The set $\mathcal{U}$ is a multiplicative
group, i.e. if $U_1(z)$ and $U_2(z)$ belong to $\mathcal{U}$
then so is the product $U_1U_2$. However, the set $\mathcal{U}$ is
not closed under addition. In contrast, the set $\mathcal{C}_N$
in \eqref{eq:m1} is closed under both addition and multiplication
by a scalar. However, $\mathcal{C}$ is not a multiplicative group,
i.e. if $F_1(z)$ and $F_2(z)$ belong to $\mathcal{C}$, the product
$F_1F_2$ does not necessarily belong to this set.
\vskip 0.2cm

The purpose of this work is to provide an easy-to-compute
characterization of all wavelet filters described by
\eqref{eq:WaveletFilters} in both presentations:
input-output and state-space. To this end we introduce the notion
of an {\em elementary wavelet filter}, i.e. a minimal McMillan
degree element in $\mathcal{W}_N$. This degree turns to be
$\frac{1}{2}N(N-1)$. A wavelet filter is then characterized as an
elementary filter multiplied from the left by an element of
$\mathcal{U}_N$ of the form $U(z^N)$ (in signal processing
terminology, \mbox{$N$-decimated} \mbox{$N$-expanded} $U$, see
e.g. \cite[Subsection 4.1.1]{Va}). In turn, this $U(z^N)$ can be
factorized into $m$ ~{\em elementary} members of $\mathcal{U}$.
We call this $m$ the {\em index}~ of the filter. The McMillan
degree of the resulting wavelet filter is
$N\left(\frac{1}{2}(N-1)+m\right)$.

\subsection{Background motivation}

As already pointed out, from application point of view wavelet
filters are of interest \cite{AJL}, \cite{AJLM}, \cite{BJ2} and
the references cited there. In their simplest form, for
discrete-time
signals, the study of filters makes use of dual frequency
variables, leading to functions of a complex variable.
\vskip 0.2cm

Aside from applications, the uses of filters in mathematics are
manifold: In duality theories in harmonic analysis, see e.g.,
\cite{HW},
\cite{JMP}, \cite{DJ1}, \cite{DJ2}, \cite{DJ3}; in wavelets
from the fundamental work of Daubechies in \cite{Da} to
\cite{Ma}, \cite{AMM}, \cite{AtGa}, \cite{SHGB}, \cite{SN},
\cite{Wo}, and other families of frames and orthogonal systems
\cite{Wa}, \cite{XM} ; in operator theory \cite{AlGo1};
%
%in operator algebras \cite{BJ2}, \cite{DJ2}, in
%representation theory (for examples the notion of filter
%bank we discuss is closely tied to certain representations
%of the Cuntz algebras), in geometry (for example what we
%introduce below as McMillan degree also arises in geometry
%as an index); the theory of matrix functions and their
%factorization, see e.g., \cite{OTHN}, \cite{WP},
%\cite{HJ1}, \cite{J1}, \cite{J2}, \cite{J3}, \cite{J4},
%\cite{J5} and the references cited there; as well as in
%mathematical physics. All of these instances are in addition
%related fields in engineering, e.g., in \cite{Be},
%\cite{BJ2}, \cite{CZZM1}, \cite{CZZM2}, \cite{CZZM3},
%\cite{CZZM4}, \cite{GVKDM} the fundamental
%ideas have been influential; as well as
%
in operator algebras \cite{BJ2}, \cite{DJ2}, in
representation theory (for examples the notion of filter
bank we discuss is closely tied to certain representations
of the Cuntz algebras), in geometry (for example what we
introduce below as McMillan degree also arises in geometry
as an index); the theory of matrix functions and their
factorization, see e.g., \cite{OTHN}, \cite{WP}, \cite{HJ1},
\cite{J1}, \cite{J2}, \cite{J3}, \cite{J4}, \cite{J5} and
the references cited there; as well as in mathematical
physics. All of these instances are in addition related
fields in engineering, e.g., in \cite{Be}, \cite{BJ2},
\cite{CZZM1}, \cite{CZZM2}, \cite{CZZM3}, \cite{CZZM4},
\cite{GVKDM} the fundamental ideas have been influential; as
well as in neighboring fields; see for example \cite{AJL},
\cite{AJLM}, \cite{BJ0}, \cite{BJ1}, \cite{TV}, \cite{Va}.
Even within harmonic analysis, there are several
viewpoints. Some authors, for example, (see e.g., \cite{CMS}
study the geometry of inner functions with the use of the
Cuntz relations built from composition operators. By
contrast, our present work takes as its point of departure
de Branges spaces and applications to systems theory, see
e.g., \cite{Ka}, \cite{AlGo1};
and control theory, see \cite{So}, \cite{VHEK}.
\vskip 0.2cm

As already mentioned, this paper addresses different audiences:
Engineers from the fields of system and control or from signal
processing and mathematician interested in classical analysis.
Thus, some of this work is of background nature (e.g. sections
\ref{sec:3}, \ref{sec:5}, \ref{sec:11}). Each section covers a
theme and begins with a brief summary of the main ideas.
\vskip 0.2cm

The outline is as follows.  In Section
\ref{sec:2} we address the set $\mathcal{C}$ described in
\eqref{eq:m1}, while in Section \ref{sec:3}, the set $\mathcal{U}$
\eqref{eq:unit}. In Section \ref{sec:4} wavelet filters are
characterized in the input-output framework. The rest of the
work is devoted to state-space realization. Relevant known
background is given for general systems in Section \ref{sec:5}
and for the set $\mathcal{U}$ in subsection \ref{Ap3}.
Realizations of elementary: Wavelet filters and unitary
functions are
given in Sections \ref{sec:6} and \ref{sec:7} respectively. The
combination of both yields the realization of wavelet filters
in \ref{sec:8}. The powerful parametrization of wavelet
filters introduced here, opens the door for future research on
both sides Engineering and Mathematics. Sample future research
topics are given in Section \ref{sec:10}. Additional background
is relegated to the appendix in section \ref{sec:11}.
\vskip 0.2cm

Many of the ideas here are familiar to some research communities.
Nevertheless, as far as we know, none of the results below,
designated as original, appeared before.

\section{The family $\mathcal{C}_N$}\label{sec:2}
\setcounter{equation}{0}

Here we outline (see Lemma \ref{Cproperty}) a practical and
geometric characterization of the set $\mathcal{C}$ defined
through the symmetry condition we introduced in \eqref{eq:m1}.

To this end, we denote by $\hat{P}$ the following $N\times N$
permutation matrix,
\[
\hat{P}=\begin{pmatrix} 0_{1\times (N-1)}&1\\
I_{N-1}& 0_{(N-1)\times 1}\end{pmatrix}.
\]
Using \eqref{eq:e} note that,
\begin{equation}\label{eq:P}
\hat{P}=U{\rm diag}\{1,~\e,~\e^2,~\ldots~,~\e^{N-1}\}U^*,
\end{equation}
for some unitary $U$, i.e. $U^*U=I$. Thus, in particular
\mbox{$\hat{P}^k_{|_{k\in[0,~N-1]}}\not=I$} and
\mbox{$\hat{P}^k_{|_{k=N}}=I$.} In addition
$\hat{P}^*\hat{P}=I$.

\begin{La}\label{Cproperty}
I. A rational function $F(z)$ is in
$\mathcal{C}_N$ if and only if it satisfies
\begin{equation}\label{eq:sym1}
F(\e z)=F(z)\hat{P}.
\end{equation}
II. Let $F_a(z)$ and $F_b(z)$ be in $\mathcal{C}_N$ then,
\[
F_b(\e z)F_a(\e z)^*=F_b(z)F_a(z)^*.
\]
If in addition $F_a(z)$ is invertible then,
\begin{equation}\label{eq:FaFb}
F_b(\e z)F_a(\e z)^{-1}=F_b(z)F_a(z)^{-1}.
\end{equation}
\end{La}

{\bf Proof:} I.\quad
One direction is clear. If
\begin{equation}\label{opera_bastille}
F(z)=\begin{pmatrix}f_0(z)&f_1(z)&
%f_3(z)&
\ldots&f_{N-1}(z)\end{pmatrix}.
\end{equation}
%$F(z)$ 
is in $\mathcal{C}_N$, that is, it is of the form \eqref{eq:m1},
then the columns of $F(\epsilon z)$ are cyclically shifted to
the left by one, namely
$f_j(z):=\left(
\begin{smallmatrix}\hat{f}_0({\epsilon}^jz)\\
\vdots\\ \hat{f}_{N-1}({\epsilon}^jz)\end{smallmatrix}\right)$
with $j=0,~\ldots~,~N-1$, so indeed \eqref{eq:sym1} holds.
%the last one becoming the first one. 
%So, $F(\epsilon z)=F(z)\widehat{P}$.
\vskip 0.2cm

Conversely, if \eqref{eq:sym1} holds, $~f_0(z),~\ldots~,~f_{N-1}(z)$
the columns of $F(z)$, see \eqref{opera_bastille}, satisfy
\[
\begin{matrix}
f_0(\epsilon z)&=&f_{1}(z)\\
f_1(\epsilon z)&=&f_{2}(z)\\
\vdots\\\
f_{N-2}(\epsilon z)&=&f_{N-1}(z)\\
f_{N-1}(\epsilon z)&=&f_0(z)
\end{matrix}
\]
and so 
$f_2(z)=f_0(\epsilon^2z),~\cdots~,~f_{N-1}(z)=f_0(\epsilon^{N-1}z)$.
Namely, $F(z)$ is of the form \eqref{eq:m1}, i.e. it belongs to
$\mathcal{C}_N$.
%
%{\bf edit}
%
%Let $F(z)$ be a $N\times N$-valued rational
%function whose columns are $f_1(z),\ldots f_N(z)$, i.e.
%\begin{equation}
%\label{opera_bastille}
%F(z)=\begin{pmatrix}f_1(z)&f_2(z)&f_3(z)&
%\ldots&f_N(z)\end{pmatrix}.
%\end{equation}
%Multiplying $F(z)$ from the right by $\hat{P}$ shifts
%the columns cyclically to the left, namely
%\[
%F(z)\hat{P}=\begin{pmatrix}f_2(z)&f_3(z)&\cdots
%&f_{N}(z)&f_1(z)\end{pmatrix}.
%\]
%On the other hand,
%\[
%F(\e z)=\begin{pmatrix}f_1(\e z)&f_2(\e z)&\ldots&f_N(\e z)
%\end{pmatrix}.
%\]
%Now, recall from \eqref{eq:m1} that each column of
%$F(z)$ in \eqref{opera_bastille} is of the form,
%\begin{center}
%$f_j(z):=
%\begin{pmatrix}
%\hat{f}_0(\e^{j-1}z)\\
%\hat{f}_1(\e^{j-1}z)\\
%\vdots\\
%\hat{f}_{N-1}(\e^{j-1}z)
%\end{pmatrix}\quad\quad\quad j=1,~\cdots~,~N.$
%\end{center}
%Substituting in Equation \eqref{eq:sym1} leads to
%\[
%\begin{split}
%\begin{pmatrix}f_1(\e z)&f_2(\e z)&\cdots
%&f_{N-1}(\e z)&f_N(\e z) \end{pmatrix} &=\\
%&\hspace{-3.cm}=
%\begin{pmatrix}f_2(z)&f_3(z)&\cdots
%&f_{N}(z)&f_1(z)\end{pmatrix}.
%\end{split}
%\]
\vskip 0.2cm

II. Substituting in \eqref{eq:sym1} and using \eqref{eq:P}
yields
\[
F_b(\e z)F_a(\e z)^*=F_b(z)\hat{P}\left(F_a(z)\hat{P}\right)^*=
F_b(z)\hat{P}\hat{P}^*F_a(z)^*=F_b(z)F_a(z)^*,
\]
so the first part is established. For the second part, use
\eqref{eq:sym1} to write \mbox{$\hat{P}=F_a(z)^{-1}F_a(\e z)$.}
Substituting in \eqref{eq:sym1} for $F_b(z)$ yields,
\[
F_b(\e z)=F_b(z)\hat{P}_{|_{\hat{P}=F_a(z)^{-1}F_a(\e z)}}=
F_b(z)F_a(z)^{-1}F_a(\e z).
\]
By assumption $F_a(\e z)$ is invertible. Multiplying by
$F_a(\e z)^{-1}$ from the right completes the proof.
\mbox{}\qed\mbox{}\\

Note that the relation in \eqref{eq:FaFb} holds whenever
\[
F_b(z)=G(z)F_a(z),
\]
for arbitrary invertible matrix valued rational function
$G(z)$ satisfying \mbox{$G(\e z)=G(z)$.} Note also that this
is the case if one takes $G$ to be $G(z^N)$ (the passage
from $G(z)$ to $G(z^N)$ is addressed in part \ref{Ap1}). It
turns out that this sufficient condition is also necessary.

\begin{Pn}\label{Pr:MR}
Let $F_a(z)$ and $F_b(z)$ be in $\mathcal{C}_N$ with
$F_a$ invertible. Then,

I. There exists a matrix valued rational function $G$ such
that
\begin{equation}
F_b(z)F_a(z)^{-1}=G(z^N).
\end{equation}
II. If in addition $G\hat{P}=\hat{P}G$ then also
\[
F_a(z)^{-1}F_b(z)=G(z^N).
\]
\end{Pn}

{\bf Proof:} I. Denoting $F:=F_bF_a^{-1}$, from
\eqref{eq:FaFb} it follows that
\begin{equation}\label{bastille1}
F(\e z)=F(z).
\end{equation}
The function $F$ has a Laurent expansion
\[
F(z)=\sum_{k=-k_0}^\infty F_kz^k
\]
converging in a punctured disk $0<|z|<r$ for some $r>0$.
Equation \eqref{bastille1}
implies that
\[
\sum_{k=-k_0}^\infty F_kz^k=\sum_{k=-k_0}^\infty F_k\e^kz^k.
\]
By uniqueness of the Laurent expansion we get that
\[
F_k=0,\quad {\rm for}\quad k\not\in N\mathbb Z.
\]
Thus, if $k_0>0$, we may assume without loss of generality that
\mbox{$k_0=Nn_0$} for some $n_0\in\mathbb N$. Thus, one
can write,
\[
F(z)=G(z^N),
\]
so part I is established.
\vskip 0.2cm

For part II denote $\hat{F}=F_a^{-1}F_b$ and note that
\[
\begin{matrix}
\hat{F}(\e z)&=&F_a(\e z)^{-1}F_b(\e z)&=&
\left(F_a(z)\hat{P}\right)^{-1}F_b(z)\hat{P}&~\\~&
=&\hat{P}^*\hat{F}(z)\hat{P}{_{|_{\hat{F}\hat{P}=
\hat{P}\hat{F}}}}&=&\hat{P}^*\hat{P}\hat{F}(z)&=
\hat{F}(z).
\end{matrix}
\]
In a way similar to part I one obtains $\hat{F}=\hat{G}(z^N).$
\mbox{}\qed\mbox{}\\

The proof of Proposition \ref{Pr:MR} can be mimicked to
obtain the following result:

\begin{Pn}
Let $F_a(z)$ and $F_b(z)$ be in $\mathcal{C}_N$. Then the
products
\[
F_b(z)\left(F_a(z^*)\right)^*\quad{\rm and}\quad
F_b(z)\left(F_a(1/z^*)\right)^*
\]
are matrix valued rational functions of $z^N$.
\end{Pn}

{\bf Proof:}
Here we make use of \eqref{eq:P} to see that for
\[
%\mbox{$
G(z)=F_b(z)\left(F_a(z^*)\right)^*
%$} 
\quad\quad{\rm or}\quad\quad
%\mbox{$
G(z)=F_b(z)\left(F_a(1/z^*)\right)^*
%$}
\]
we have
\eqref{bastille1} i.e. \mbox{$G(z)=G(\e z)$.}
\mbox{}\qed\mbox{}\\

We now introduce the notion of elementary wavelet filters. To
this end we need some preliminaries. For $N\geq 2$ let
$Q_N$ be the following (constant) $N\times N$ unitary matrix,
\begin{equation}\label{eq:Fourier}
Q_N:=\frac{1}{\sqrt{N}}
\begin{pmatrix}
\e^{-(0\cdot 0)}&\e^{-(0\cdot 1)}&\e^{-(0\cdot 2)}&\cdots
&e^{-0\cdot(N-1)}\\
\e^{-(1\cdot 0)}&\e^{-(1\cdot 1)}&\e^{-(1\cdot 2)}&
\cdots&e^{-1\cdot(N-1)}\\
\e^{-(2\cdot 0)}&\e^{-(2\cdot 1)}&\e^{-(2\cdot 2)}&
\cdots&e^{-2\cdot(N-1)}\\
\vdots& \vdots& \vdots& \vdots&\vdots\\
\e^{-(N-1)\cdot 0}&\e^{-(N-1)\cdot 1}&\e^{-(N-1)\cdot 2)}&
\cdots&e^{-(N-1)\cdot(N-1)}
\end{pmatrix}.
\end{equation}
For example,
\begin{center}
$Q_{|_{N=4}}=${\mbox{\tiny$\frac{1}{2}$}}
$\begin{pmatrix}1&~~~1&~~~1&~~~1\\ 1&-i&-1&~~i\\
1&-1&~~~1&-1\\ 1&~~~i&-1&-i\end{pmatrix}$.
\end{center}
\vskip 0.2cm

Recall that $Q_N$ is called the Discrete Fourier Transform
matrix,
see e.g. \cite[Subsection 3.1.1]{BN}, \cite[page 794]{Va}.
Namely, if $\hat{x}\in\C^N$
is the discrete Fourier transform of $x\in\C^N$, then
they are related through
\begin{center}
$\hat{x}=${\mbox{\tiny$\sqrt{N}$}}$Q_Nx\quad\quad\quad\quad
$x=${\mbox{\tiny$\frac{1}{\sqrt{N}}$}}Q_N^*\hat{x}.$
\end{center}

Now, we denote by $\hat{U}(z)$ the following special element of
$\mathcal{U}_N$,
\begin{equation}\label{eq:hatU}
\hat{U}_N(z):={\rm diag}\left\{1,~\frac{1}{z}~,~\ldots~,~
\frac{1}{z^{N-1}}\right\}.
\end{equation}
Note that $\hat{U}_N(z)$ is of McMillan degree $\frac{1}{2}N(N-1)$.
\vskip 0.2cm

Next, using \eqref{eq:Fourier} and \eqref{eq:hatU} for
$N\geq 2$, we define the {\em elementary wavelet filter,}
\begin{equation}\label{eq:ElementaryW}
\hat{W}_N(z):=\hat{U}_N(z)Q_N~.
\end{equation}
Note that the McMillan degree of $\hat{W}_N$ is as of
$\hat{U}_N$, i.e. $\frac{1}{2}N(N-1)$.

\begin{Cy}\label{Cy:FamilyC}
I. Let $F\in\mathcal{C}_N$ be given and let $L(z)$, $R(z)$ be
arbitrary $N\times N$ valued rational functions, almost
everywhere invertible. Then\begin{footnote}{Below,
$L$ and $R$ stand for ``left" and ``right" respectively.}
\end{footnote}
\[
\{L(z^N)F(z)\}\quad\quad{\rm and}\quad\quad
\{F(z)R(z^N)~:~R(z^N)\hat{P}=\hat{P}R(z^N)\},
\]
are families in $\mathcal{C}_N$.

II. For given $N\geq 2$ let $\hat{W}_N(z)$ be the elementary
wavelet filter in \eqref{eq:ElementaryW}. Then,
\[
\hat{W}\in\mathcal{C}.
\]
Moreover, all elements in $\mathcal{C}_N$ are given by
\begin{equation}\label{eq:CharC}
\mathcal{C}_N=\{~L(z^N)\hat{W}_N(z)~:~L\quad{\rm rational,~
invertible}~\}.
\end{equation}
\end{Cy}
\vskip 0.2cm

Item I follows from Proposition \ref{Pr:MR}.

II. To see that $\hat{W}$ in \eqref{eq:ElementaryW} is indeed
in $\mathcal{C}_N$, note that
\[
\hat{W}(\e z)=
\hat{U}(\e z)Q=
\hat{U}(z)\left({\hat{U}(z)}_{|_{z=\e}}\right)Q=
\hat{U}(z)QP=\hat{W}(z)\hat{P}.
\]
Now \eqref{eq:CharC} follows from item I since
%$\det\hat{W}_N(z)\not=0$.
$\hat{W}_N(z)$ is invertible.

\section{The group $\mathcal{U}$}
\label{sec:3}
\setcounter{equation}{0}

As already mentioned the infinite-dimensional group
$\mathcal{U}$ is well known in literature see e.g. \cite{AlGo1},
\cite{BJ1}, \cite{DJ2}, \cite{J5}
and the references cited there. We here review properties
of this set, necessary to our construction in the sequel.
\vskip 0.2cm

Recall that if $U(z)$ is in the set $\mathcal{U}$, then so
are $\left(U(z^*)\right)^*$ and if stability is compromised
also $U(\frac{1}{z})$.
\vskip 0.2cm

To construct the set $\mathcal{U}$ first recall that (up to
multiplication by a constant of a unity modulus) a unity
degree, scalar, stable rational function $\phi_{\alpha}(z)$
mapping $\mathbb T$ to itself can always be written as,
\begin{equation}\label{ElemntaryPhi}
\phi_{\alpha}(z)=\frac{1-{\alpha}^*z}{z-\alpha}\quad\quad
\alpha\in\mathbb{D}.
\end{equation}
We shall call $\phi_a(z)$ in \eqref{ElemntaryPhi} an {\em
elementary scalar unitary} function. Every scalar rational
function mapping $\mathbb{T}$ to itself of McMillan degree
$p$, may be factored to {\em elementary} functions of the
form \eqref{ElemntaryPhi}, namely
\[
\prod\limits_{j=1}^p\frac{1-{\alpha_j}^*z}{z-\alpha_j}\quad\quad
\alpha_j\in\mathbb{D},
\]
where the $\alpha_j$'s are not necessarily distinct. As in
this rational function the poles are within $\mathbb{D}$ and
the zeros are outside $\overline{\mathbb D}$, no cancellation
can occur. Thus the issue of minimality of realization, central
to system theory see e.g. \cite{Ka}, \cite[Section 5.5]{So},
can here be avoided.
\vskip 0.2cm

We next extend this idea to factorize elements in
$\mathcal{U}_N$. To this end, we need some preliminaries.
The unit sphere in $\C^N$ will be denoted by,
\[
\mathbb{S}_{N-1}:=\{ v\in\C^N~:~v^*v=1\}.
\]
We can now introduce an $N$-dimensional {\em elementary
unitary matrix} of McMillan degree one,
\begin{equation}\label{eq:U(u,c)}
V(z):=I_N+\left(\frac{1-{\alpha}^*z}{z-\alpha}-1\right)vv^*,
\quad\quad v\in\mathbb{S}_{N-1}\quad\alpha\in\mathbb{D}.
\end{equation}
Recall that $V(z)$ is a Finite Impulse Response (FIR) filter
if and only if \mbox{$\alpha=0$.}
\vskip 0.2cm

For $j$ natural, we shall find convenient to denote,
\begin{equation}\label{eq:Uj(u,c)}
V_j(z):=V(z, v_j, \alpha_j),
\end{equation}
with $v_j\in\mathbb{S}_{N-1}$ and ${\alpha}_j\in\mathbb{D}$,
parameters.
\vskip 0.2cm

We can cite an adapted version of \cite[Theorem 3.11]{AlGo1},
see also \cite[Eq. (14.9.2)]{Va}

\begin{Tm}\label{Tm:U(z)}
Using \eqref{eq:Uj(u,c)} (up to multiplication by constant unitary
matrices from the left and from the right) $U(z)$ in $\mathcal{U}_N$
of McMillan degree $~p$, can always be written as
\begin{equation}\label{eq:ProdU_j(z)}
U(z)=\prod\limits_{j=1}^pV_j(z).
\end{equation}
\end{Tm}

We now illustrate the richness associated with $U(z)$ in
\eqref{eq:ProdU_j(z)} already for $p=2$,
\[
U(z)=V(z, v_1, \alpha_1)V(z, v_2, \alpha_2)
%_{|_{\alpha_1=\alpha_2=\alpha}}
%U(z)=V(z^q, v_1, \alpha_1)V(z^p, v_2, \alpha_2)_{|_{\alpha_1=\alpha_2=\alpha}}
%{(z^p-{\alpha}_1)(z^q-{\alpha}_2)}
\]
Consider two ``extreme" cases:
\[
\begin{matrix}
{U(z)}_{ |_{v_1=v_2=v}}&=&
I_N+\left(\frac{(1-{\alpha}_1^*z)(1-{\alpha}_2^*z)}
{(z-{\alpha}_1)(z-{\alpha}_2)}
-1\right)vv^*_{|_{\alpha_1=-\alpha_2=\sqrt{\alpha}}}\\~&
=&I_N+\left(\frac{1-{\alpha}^*z^2}{z^2-\alpha}-1\right)vv^*
\\~&=&V(z^2,~v,~\alpha),
\end{matrix}
\]
see \eqref{eq:tildeU} below. In contrast,
\begin{equation}\label{eq:uu}
\begin{matrix}
{U(z)}_{ |_{\begin{smallmatrix}v_1^*v_2=0\\
\alpha_1=\alpha_2=\alpha\end{smallmatrix}}}
%&=&I_N+
%\left(\frac{1-{\alpha}^*z}{z-{\alpha}}-1\right)v_1v_1^*
%+{\left(\frac{1-{\alpha}^*z}
%{z-{\alpha}}-1\right)v_2v_2^*}_{|_{p=q}}
%\\~
&=&I_N+\left(\frac{1-{\alpha}^*z}
{z-\alpha}-1\right){(v_1v_1^*+v_2v_2^*)}_{|_{N=2}}
\\~&=&\frac{1-{\alpha}^*z}{z-\alpha}I_2~.
\end{matrix}
\end{equation}
\vskip 0.2cm

In particular one can use Theorem \ref{Tm:U(z)} to
construct $\hat{U}(z)$ in \eqref{eq:hatU} by taking in
\eqref{eq:uu} the parameters $\alpha_j=0$ and $v_j=e_j$, the
unit vectors with the appropriate multiplicity.
Recall that the resulting $\hat{U}_N(z)$ is of McMillan degree
$\frac{1}{2}N(N-1)$ $(=$the number of elementary unitary
matrices involved).
\vskip 0.2cm

Before proceeding, we find it convenient to introduce
the $\rho$-scaled disk
\[
\mathbb{D}_{\rho}:=\rho\cdot\mathbb{D}
\quad\quad\quad\quad\rho\in[0,~1].
\]
We shall let $\alpha$, the pole of $V(z)$ in \eqref{eq:U(u,c)} 
to lie within $\mathbb{D}_{\rho}$ and thus to quantify the
stability of the unitary matrices at hand: From arbitrary
asymptotic (in fact exponential) stability corresponding to
$\rho=1$, to finite impulse response, when $\rho=0$. We shall
thus call $\rho$ ~{\em the spectral radius}~ of the unitary
matrix $V(z)$. The spectral radius, $\rho$, is associated with
concept of spectrum of an operator, for further motivation
and details see \eqref{eq:SpectRad}. A word of caution, this
is not to be confused with notion of the frequency spectrum of
signals, commonly used in electrical engineering.
\vskip 0.2cm

Following \eqref{eq:U(u,c)} and \eqref{eq:ProdU_j(z)} we find it
convenient to introduce an \mbox{$N\times N$}, {\em
$N$-decimated $N$-expanded elementary unitary}~ matrix of McMillan
degree $N$,
\begin{equation}\label{eq:tildeU}
V(z^N):=I_N+\left(\frac{1-{\alpha}^*z^{N}}{z^{N}-{\alpha}}-1\right)vv^*,
\end{equation}
with $v_j\in\mathbb{S}_{N-1}$ and ${\alpha}_j\in\mathbb{D}_{\rho}$
parameters, for some prescribed $\rho\in[0,~1].$
\vskip 0.2cm

For future reference we now formulate a version of Theorem
\ref{Tm:U(z)} for $N$-decimated $N$-expanded unitary matrix.

\begin{Cy}\label{Cy:U(zN)}
Let $U\in\mathcal{U}_N$ be of the form $U(z^N)$. Then
(up to multiplication by constant unitary matrices from
the left and from the right) using \eqref{eq:tildeU} it can
always be written as,
\begin{equation}\label{eq:U(zN)}
U(z^N)=\prod\limits_{j=1}^mV_j(z^N).
\end{equation}
We shall call $m$ the index of $U(z^N)$.
Moreover, $U(z^N)$ is of McMillan degree $mN$.
\end{Cy}

With a slight abuse of terminology we shall refer to
$V(z^N)$ in \eqref{eq:tildeU} as an {\em elementary
unitary}~ matrix (omitting the ``$N$-decimated
$N$-expanded" qualifier).
\vskip 0.2cm

In the next section we combine the above properties of
the set $\mathcal{U}$ along with the set $\mathcal{C}$
from Section \ref{sec:2}, to characterize $\mathcal{W}$,
the set of wavelet filters.

\section{The set of Wavelet filters $\mathcal{W}_N$}
\setcounter{equation}{0}
\label{sec:4}

Our first observation on the structure of $\mathcal W_{N}$, the
set of $N\times N$ wavelet filters \eqref{eq:WaveletFilters}, is an
immediate consequence of Proposition \ref{Pr:MR}.

\begin{Pn}\label{Pr:WF}
Let $W_a(z)$ and $W_b(z)$ be in $\mathcal W_{N}$. Then there
exists \mbox{$U\in\mathcal{U}_N$} such that
\[
W_b(z)=U(z^N)W_a(z).
\]
\end{Pn}
\vskip 0.2cm

Indeed, substitute unitary matrices in Proposition \ref{Pr:MR}
so that, \mbox{$F_a(z)=W_a(z)$,} $F_b(z)=W_b(z)$, and
$G(z^N)=U(z^N)$.
\vskip 0.2cm

We next exploit Proposition \ref{Pr:WF} and Corollaries
\ref{Cy:FamilyC}, \ref{Cy:U(zN)} to provide an easy-to-compute
characterization of all wavelet filters in $\mathcal{W}_N$
($N$ prescribed) taking the McMillan degree as a parameter.
Strictly speaking the parameter will be $m$, the index of
the filter, yet to be defined. The exact relation between
$m$ and the McMillan degree is given in Corollary
\ref{Cy:degree} below.
\vskip 0.2cm

One can refine this parametrization and for a prescribed
$~m$, to limit the spectral radius of the wavelet
filter to $\rho$ for some $\rho\in[0,~1]$.

\begin{Tm}
I. Let $W\in\mathcal{W}_N$ be a given wavelet filter, and
let $L(z^N)$ and $R(z^N)$ be in $\mathcal{U}_N$. Then,
\[
\{L(z^N)W(z)\}\quad\quad{\rm and}\quad\quad
\{W(z)R(z^N)~:~R(z^N)\hat{P}=\hat{P}R(z^N)\},
\]
are families of wavelet filters.

II. For given $N\geq 2$ let $\hat{W}_N(z)$ be the
elementary Wavelet filter in \eqref{eq:ElementaryW}. Then,
\[
\hat{W}_N\in\mathcal{W}_N.
\]
III. All wavelet filters are
given by, \begin{equation}\label{CharHatW}
\mathcal{W}_N=\left\{\left(\prod\limits_{j=1}^m
V_j(z^N)\right)\hat{W}_N(z)~:~v_j\in\mathbb{S}_{N-1},~
\alpha_j\in\mathbb{D}_{\rho}
\right\},
\end{equation}
where $\rho\in[0,~1]$, the spectral radius of the filter, is
prescribed.

We shall call $~m~$ in \eqref{CharHatW}, the index of the filter.
\end{Tm}
\vskip 0.2cm

For an arbitrary given $W(z)$ in $\mathcal{W}_N$ we consider
equivalent all wavelet filters of the form
\[
QW(z),
\]
where $Q$ varies over all ~{constant}~ $N\times N$ unitary matrices.
\vskip 0.2cm

From \eqref{CharHatW} we have the following.
\begin{Cy}
Let $W_1(z), W_2(z), W_3(z)$ be in $\mathcal{W}_N$, then
\[
W_1W_2^{-1}W_3\in\mathcal{W}_N~.
\]
\end{Cy}
\vskip 0.2cm
The characterization of wavelet filters in \eqref{CharHatW}
implies that these are matrix valued functions of ``quantized"
McMillan degrees.

\begin{Cy}\label{Cy:degree}
An index $m$ wavelet filter of dimension $N$ is of
McMillan degree,~~$N\left(\frac{1}{2}(N-1)+m\right).$
\end{Cy}

Indeed, as already mentioned, the McMillan degree of
$\hat{W}_N(z)$ is $\frac{1}{2}N(N-1)$. From the construction in
\eqref{eq:tildeU} and \eqref{CharHatW}, each of the $m$
elementary unitary matrices $V_j(z^N)$
contributes $N$ to the overall McMillan degree.
\vskip 0.2cm

We now wish to ``count" how many wavelet filters are there.

\begin{Ob}\label{Ob:parameters}
I. An index $m$ wavelet filter of dimension $N$ and spectral
radius $\rho$,
can always be parameterized as a point in the real set,
\[
\left([0,~\pi)\times[0,~2\pi)^{2(N-1)}\times[0,~\rho)\right)^m
\quad\quad\quad\rho\in(0,~1].
\]
The set of Finite Impulse Response (FIR) filters
(=zero spectral radius) may be parameterized by,
\[
\left([0,~\pi)\times[0,~2\pi)^{2N-3}\right)^m.
\]
II. In the above parametrization, the set of wavelet filters
is convex.
\end{Ob}

{\bf Proof}\quad I.~
First, since the structure of
$\hat{W}_N$ is fixed, from \eqref{eq:ElementaryW} it is clear
that the number of free parameters in $W$ is given by the
the number of parameters in a unitary matrix $U(z^N)$ in
\eqref{eq:U(zN)}. In turn, recall that each of the $m$ elementary
unitary matrices $V_j$ in \eqref{eq:U(zN)} may be identified
with a point in $\mathbb{S}_{N-1}\times[0,~\rho)$. Thus, the
number of parameters in a unitary matrix $U(z^N)$ is
given by
\[
\left(\mathbb{S}_{N-1}\times[0,~\rho)\right)^m,
\]
and out of them, the Finite Impulse Response (FIR) filters
are given by the factor,
\[
\mathbb{S}_{N-1}^m~.
\]
Recall now that a point on $~v\in\mathbb{S}_{N-1}$, can be
equivalently parameterized as a point in the real ``box"
$[0, 2\pi)^{2N-1}$. For example,
for $N=3$ one has $v=${\mbox{\tiny$\begin{pmatrix}
\cos(\delta)e^{i\alpha}\\
\cos(\eta)\sin(\delta)e^{i\beta}\\
\sin(\eta)\sin(\delta)e^{i\gamma}\end{pmatrix}$}} with
$\alpha, \beta, \gamma, \delta, \eta\in[0, 2\pi)$.
\vskip 0.2cm

As we are actually interested only in elements of the form $vv^*$,
the number of parameters is reduced to
$[0,~\pi)\times[0, 2\pi)^{2N-3}$. For example, in the above
case of $N=3$, without loss of generality one can take
$\alpha=0$ and $\delta\in[0,~\pi)$.
\vskip 0.2cm

To complete the construction use the polar presentation
of elements within a disk ${\mathbb D}_{\rho}$ in $\C$ as
points in $[0,~2\pi)\times[0,~\rho).$
Thus, an index 1 filter is parameterized by
\[
[0,~\pi)\times[0, 2\pi)^{2N-3}\times[0,~2\pi)\times[0,~\rho)
%\frac{1}{2}\times[0, 2\pi)^{2(N-1)}\times[0,~2\pi)\times[0,~\rho]
=[0,~\pi)\times[0, 2\pi)^{2(N-1)}\times[0,~\rho)
\]
so the construction is complete.
\vskip 0.2cm

II. Is an immediate consequence of part I.
\qed
\vskip 0.2cm

The favorable properties of FIR  filters are well known,
see e.g. \cite[Sections 2.4.2, beginning of 3.3]{Va}.
In particular, they dramatically attenuate any noise.
However, Observation \ref{Ob:parameters} suggests that whenever
the noise level is sufficiently low, it is not recommended to
restrict the discussion to FIR wavelet filters: For
\mbox{$1>>\rho>0$} the filter is ``almost FIR" but there are
already ``many more" filters than FIR's. For a related
discussion, see also the end of Section III in \cite{HHN}.
\vskip 0.2cm

The convex parametrization of all wavelet filters introduced in
Observation \ref{Ob:parameters} may be in particular convenient
for optimization, see e.g. \cite{CZZM4}, \cite{HHN}, \cite{TV}
and \cite{VHEK}.
\vskip 0.2cm

Following \eqref{CharHatW} by multiplying from the left a given
$W_a\in\mathcal{W}_N$ of McMillan degree
$N\left(\frac{1}{2}(N-1)+m\right)$,
by an elementary unitary matrix $V(z^N)$, i.e. increasing the
index by 1, one can construct another wavelet filter $W_b(z)$
of McMillan degree \mbox{$N(\frac{1}{2}(N+1)+m)$.} Indeed,
\[
W_b(z)=V(z^N)W_a(z),
\]
where $V$ is as in \eqref{eq:tildeU} with
$v\in\mathbb{S}_{N-1}$ and $\alpha\in\mathbb{D}_{\rho}$,
parameters.

This is illustrated next.

\begin{Ex}
{\rm
Following \eqref{eq:ElementaryW} let $\hat{W}_2(z)$ be the
two dimensional elementary wavelet filter,
\begin{equation}\label{eq:hatW2}
\hat{W}_2(z)=\frac{1}{\sqrt{2}}\begin{pmatrix}1&~~1\\
\frac{1}{z}&-\frac{1}{z}\end{pmatrix}.
\end{equation}
Let $V_{\alpha}(z^2)$ in $\mathcal{U}_2$ be an
elementary unitary matrix of the form of \eqref{eq:tildeU}
with $v=${\mbox{\tiny$\begin{pmatrix}0\\ 1\end{pmatrix}$}}
and a parameter $\alpha\in\mathbb{D}_{\rho}$,
\begin{equation}\label{eq:Ua}
V_{\alpha}(z^2)=\begin{pmatrix}1&0\\ 0&
\frac{1-{\alpha}^*z^2}{z^2-\alpha}
\end{pmatrix}.
\end{equation}
Clearly, it is of McMillan degree 2.
\vskip 0.2cm

Multiplying $\hat{W}_2(z)$ from the left by $U_{\alpha}(z^2)$
yields $W_a(z)$, a wavelet filter with $N=2$ and $m=1$, i.e.
of McMillan degree three,
\begin{equation}\label{eq:deg3filter}
W_a(z)=V_{\alpha}(z^2)\hat{W}_2(z)=\frac{1}{\sqrt{2}}
\begin{pmatrix}1&~~1\\
\frac{1-{\alpha}^*z^2}{z(z^2-\alpha)}&
\frac{{\alpha}^*z^2-1}{z(z^2-\alpha)}
\end{pmatrix}.
\end{equation}
Next, let $V_{\beta}(z^2)$ be a member of $\mathcal{U}_2$
of index two (i.e. a product of a pair of elementary matrices of
the form of \eqref{eq:tildeU} both with $v=${\mbox{\tiny$
\begin{pmatrix}1\\ 0\end{pmatrix}$}}) where
$\beta\in\mathbb{D}_{\rho}$ is a parameter,
\begin{equation}\label{eq:Ub}
V_{\beta}(z^2)=\begin{pmatrix}\frac{1-{\beta}^*z^4}{
z^4-\beta}&0\\ 0&1\end{pmatrix}.
\end{equation}
Clearly, it is of McMillan degree 4.
\vskip 0.2cm

In turn, multiplying $W_a(z)$ in \eqref{eq:deg3filter} from the
left by $U_{\beta}(z^2)$ from \eqref{eq:Ub}, yields $W_b(z)$, a
wavelet filter with $N=2$ and $m=3$, i.e. of McMillan degree
seven (see Corollary \ref{Cy:degree}),
\begin{equation}\label{eq:deg7filter}
W_b(z)=V_b(z^2)W_a(z)=\frac{1}{\sqrt{2}}
\begin{pmatrix}\frac{1-{\beta}^*z^4}{z^4-\beta}&
\frac{1-{\beta}^*z^4}{z^4-\beta} \\~&~\\
\frac{1-{\alpha}^*z^2}{z(z^2-\alpha)}&\frac{{\alpha}^*z^2-1
}{z(z^2-\alpha)}\end{pmatrix}.
\end{equation}
Note that it is only for $\alpha=\beta=0$ that $W_b(z)$ is a
FIR filter, see
%Moreover it is only for $1>|\alpha |,|\beta |$ that this filter
%is Schur stable, see
Observation \ref{Ob:parameters}.
\vskip 0.2cm

State space realization of $\hat{W}_2$ in equation
\eqref{eq:hatW2}, is given in item 1 of Example
\ref{Ex:canonicalM} below. State space realizations of
$W_a(z)$ and $W_b(z)$ in equations  \eqref{eq:deg3filter} and
\eqref{eq:deg7filter}, respectively are given in Example
\ref{Ex:RealizationWaWb} below.
}\qed
\end{Ex}

The rest of this work is devoted to state-space realization of
wavelet filters.  This topic has gained popularity in the last
two decades, see e.g.  \cite{CZZM1}, \cite{CZZM2}, \cite{CZZM3},
\cite{CZZM4}, \cite{GVKDM}, \cite{SHGB}, \cite{TV},
\cite[Section 13.4]{Va}, \cite{VHEK} and \cite{Wa}.
In Sections \ref{sec:5}, \ref{sec:6} below we review the
background necessary to translate the above presented
construction procedure of wavelet filters from matrix
valued rational functions setting (input-output in
engineering terminology) to state-space framework.

\section{State space realization - preliminaries}
\label{sec:5} \setcounter{equation}{0}

Every $N\times N$ valued rational function $F(z)$ analytic
at infinity ($\lim\limits_{z\rightarrow\infty}F(z)$
exists) can be written as
\[
F(z)=C(zI_p-A)^{-1}B+D,
\]
where $B\in\C^{p\times N}$ and $C\in\C^{N\times p}$. This may
be viewed as the $Z$-transform of
\begin{equation}\label{realization1}
x(n+1)=Ax(n)+Bu(n)\quad\quad\quad y(n)=Cx(n)+Du(n),
\end{equation}
where $u, y$ are $N$-dimensional input and output
respectively and $x$ is $p$-dimensional (state in system
theory terminology). Namely, $F(z)$
maps the $Z$-transform of $u$ to the $Z$-transform of $y$.
We find it
convenient to write \eqref{realization1} in the
\mbox{$(p+N)\times(p+N)$} {\sl system matrix} form,
\begin{equation}\label{eq:SystemMatrix}
M=\left(\begin{array}{r|r}A&B\\ \hline C&D\end{array}\right).
%{smallmatrix}
\end{equation}
In principle, rectangular system matrices, (i.e. the number
of inputs differs from the number of outputs) are of interest,
but they are mostly beyond the scope of this work,
see Remark 6 in Section \ref{sec:10}. As already mentioned,
the issue of minimality of realization is here trivially
satisfied and thus can be avoided.
\vskip 0.2cm

Before proceeding, recall that $\rho(A)$, the spectral radius
of a square matrix $A$, is the radius of the smallest disk,
centered at the origin of the complex plane, containing the
eigenvalues (=spectrum) of $A$. For some details, see e.g.
\cite{HJ1} items 5.6.8-5.6.14. In particular, recall that for
in \eqref{realization1} one has that 
there exists $\beta\geq 1$ so that for all natural $n$,
\begin{equation}\label{eq:SpectRad}
{\| x(n)\|}_{|_{u(n)\equiv 0}}\leq\beta\| x(0)\|{\rho}^n,
\quad\quad\quad\quad\forall x(0)
\end{equation}
where $\rho=\rho(A)$, the spectral radius of the matrix $A$,
which in turn coincides with the spectral radius defined in
Sections \ref{sec:3}, \ref{sec:4}.
\vskip 0.2cm

For future reference we now recall in the state space realization
of a composition of a pair of rational functions
(series connection or cascade of systems in Electrical Engineering
terminology). Namely $F_{\Delta}(z)$, $F_a(z)$ are of compatible
dimensions and $F_b(z)$ is obtained,
\begin{equation}\label{eq:Fb}
F_b(z)=F_{\Delta}(z)F_a(z).
\end{equation}
Assuming the state-space realization of each $F_a(z)$ and
$F_{\Delta}(z)$ is known, one can construct a realization of
the resulting $F_b(z)$, see e.g. \cite[Subsection 8.3.3]{Ka},
\cite[Eq. (4.15)]{Wa}.

\begin{Ob}\label{Ob:series}
Given $l\times q$ and $q\times r$ valued rational functions
$F_{\Delta}(z)$, $F_a(z)$, respectively, admitting
state space realization
\[
F_a(z)=C_a(zI-A_a)^{-1}B_a+D_a\quad\quad\quad
A_a\in\C^{p_a\times p_a}
\]
\[
\begin{matrix}
B_a\in\C^{p_a\times r}&~&~&
C_a\in\C^{q\times p_a}&~&~&
D_a\in\C^{q\times r}
\end{matrix}
\]
\[
F_{\Delta}(z)=C_{\Delta}(zI-A_{\Delta})^{-1}B_{\Delta}+D_{\Delta}
\quad\quad\quad
A_{\Delta}\in\C^{p_{\Delta}\times p_{\Delta}}
\]
\[
\begin{matrix}
B_{\Delta}\in\C^{p_{\Delta}\times q}&~&~&
C_{\Delta}\in\C^{l\times p_{\Delta}}&~&~&
D_{\Delta}\in\C^{l\times q}.
\end{matrix}
\]
The system matrix $M_b$ associated with the realization
of $F_b(z)$ in \eqref{eq:Fb} is given by,
\[
M_b=\left(
\begin{array}{rr|r}
%\begin{smallmatrix}
A_{\Delta}&B_{\Delta}C_a
%&\vdots
&B_{\Delta}D_a\\0&~A_a
%&\vdots
&B_a\\
%\cdots&\cdots&\cdots&\cdots\\
\hline
C_{\Delta}&~D_{\Delta}C_a
%&\vdots
&D_{\Delta}D_a
\end
%{smallmatrix}
{array}\right).
\]
\end{Ob}
\vskip 0.2cm

Although not central to our discussion, we comment that even when
the realizations of $F_a(z)$ and $F_{\Delta}(z)$ are minimal, the
realization in $M_b$ is not necessarily minimal.
\vskip 0.2cm

We now point out that whenever $A_a$ and $A_{\Delta}$ are upper
triangular, the resulting $A_b$ is upper triangular as well.
Fortunately, this is indeed the case with all wavelet filters,
see e.g. \eqref{eq:Mb}, the system matrix associated with
$W_b(z)$ from
\eqref{eq:deg7filter}. See also Remark 2 in Section \ref{sec:10}.
\vskip 0.2cm

State space realization of unitary functions are briefly
discussed in part \ref{Ap3}.
\vskip 0.2cm

We now outline the rest of this work: From the characterization
of wavelet filters in \eqref{CharHatW} it follows that every
wavelet
filter of dimension $N$ and index $m$ is a product of the form
\[
\prod_{j=1}^mV_j(z^N)\hat{W}(z).
\]
State space realizations of elementary wavelet filters
$\hat{W}_N(z)$ and of elementary unitary matrices
$V_j(z^N)$ are given in sections \ref{sec:6} and \ref{sec:7}
respectively. In section \ref{sec:8} we combine them.
\vskip 0.2cm

\section{State space realization of elementary wavelet
filters}\label{sec:6} \setcounter{equation}{0}

Recall that $\hat{W}_N(z)$, the elementary wavelet filter
in \eqref{eq:ElementaryW} can be described as the minimal
McMillan degree $(=\frac{1}{2}N(N-1)~)$ filter in
$\mathcal{W}_N$. In addition, it is a FIR filter. Thus,
there is no surprise that there is a two stages recipe for
the construction of $\hat{M}_N$, the
associated system matrix.
\vskip 0.2cm

1. We start by constructing an auxiliary permutation matrix,
e.g. \cite[0.9.5]{HJ1}, denoted by $\tilde{M}_N$ so that
\[
\tilde{M}_N=\left(\begin{array}{r|r}
A&\tilde{B}\\ \hline C&\tilde{D}\end{array}
\right).
\]
First, $A$ is the following
$\frac{1}{2}N(N-1)\times\frac{1}{2}N(N-1)$ matrix,
\[
A={\rm diag}\{J_1(0),~\ldots~,~J_{N-1}(0)\},
\]
where $J_k(0)$ is a $k$-dimensional Jordan block
associated with a zero eigenvalue, e.g. \cite[3.1.1]{HJ1}.
Next, $\tilde{D}$ is a $N\times N$ matrix with zeroes
everywhere, except a single 1 in the upper left corner.
\vskip 0.2cm

Now, $C$ and $\tilde{B}$ are constructed so that
(i) the 1's in each are in a descending
staircase order, e.g.
\[
\left(\begin{smallmatrix}0&0&0&0&0&0\\
1&0&0&0&0&0\\
0&1&0&0&0&0\\
0&0&0&1&0&0\end{smallmatrix}\right)
\]
and (ii) the resulting $\tilde{M}_N$ is a permutation matrix.

2. Using $Q_N$ from \eqref{eq:Fourier},
\[
\hat{M}_N=\tilde{M}_N{\rm diag}\{I_p,~Q_N\}
=\left(\begin{array}{r|r}
A&\tilde{B}Q_N\\ \hline C&\tilde{D}\hat{V}_N\end{array}
\right)=\left(\begin
{array}{r|r}A&B\\ \hline C&D\end{array} \right).
\]
\vskip 0.2cm

This is illustrated by the following.

\begin{Ex}\label{Ex:canonicalM}
{\rm
We here consider the system matrix $\hat{M}_N$ \eqref{eq:SystemMatrix}
associated with the elementary wavelet filter $\hat{W}_N(z)$ in
\eqref{eq:ElementaryW}.
\vskip 0.2cm

1. For $N=2$, it is given by,
\[
\tilde{M}_N=\left(\begin{array}{r|rr}
0     & 0  &   1  \\ \hline 0     & 1   &   0  \\
1     & 0   &   0 \end{array}\right)\quad\quad\quad
\hat{M}_N=\left(\begin{array}{r|rr}
0&\frac{1}{\sqrt{2}}&-\frac{1}{\sqrt{2}}\\
\hline 0 &\frac{1}{\sqrt{2}}&~\frac{1}{\sqrt{2}}\\
1 &         0        &   ~~0\end{array}\right).
\]
\vskip 0.2cm

2. For $N=4$, the system matrix is given by,
\[
\tilde{M}_N=\left(\begin{array}{rrrrrr|rrrr}
0&0&0&0&0&0&0&1&0&0\\
0&0&1&0&0&0&0&0&0&0\\
0&0&0&0&0&0&0&0&1&0\\
0&0&0&0&1&0&0&0&0&0\\
0&0&0&0&0&1&0&0&0&0\\
0&0&0&0&0&0&0&0&0&1\\
\hline
0&0&0&0&0&0&1&0&0&0\\
1&0&0&0&0&0&0&0&0&0\\
0&1&0&0&0&0&0&0&0&0\\
0&0&0&1&0&0&0&0&0&0
\end{array}\right)\quad\quad
\hat{M}_N=\left(\begin{array}{rrrrrr|rrrr}
0&0&0&0&0&0
&\frac{1}{2}&\frac{-i}{2}&\frac{-1}{2}&\frac{i}{2} \\
0&0&1&0&0&0
&  0        &     0      &    0       &    0       \\
0&0&0&0&0&0
&\frac{1}{2}&\frac{-1}{2}&\frac{1}{2} &\frac{-1}{2}\\
0&0&0&0&1&0
&  0        &     0      &    0       &    0       \\
0&0&0&0&0&1
&  0        &     0      &    0       &    0       \\
0&0&0&0&0&0
&\frac{1}{2}&\frac{i}{2} &\frac{-1}{2}&\frac{-i}{2}\\
\hline
0&0&0&0&0&0
&\frac{1}{2}&\frac{1}{2} &\frac{1}{2} &\frac{1}{2} \\
1&0&0&0&0&0
&  0        &     0      &    0       &    0       \\
0&1&0&0&0&0
&  0        &     0      &    0       &    0       \\
0&0&0&1&0&0
&  0        &     0      &    0       &    0
\end{array}\right).
\]
}\qed\end{Ex}

\section{State space realization of elementary unitary matrices}
\label{sec:7} \setcounter{equation}{0}

We start by showing that the system matrix associated with an
elementary unitary {\em matrix-valued} function see \eqref{eq:tildeU},
may be obtained from the system matrix associated with an elementary
{\em scalar} function.

\begin{Ob}\label{Ob:MtildeU}
Let $V(z^N)$, be an elementary unitary function as in
\eqref{eq:tildeU}, i.e.
\begin{equation}\label{eq:u(z^N)}
V(z^N)=
I_N+\left(\frac{1-{\alpha}^*z^N}{z^N-\alpha}-1\right)vv^*
\quad\quad v\in\mathbb{S}_{N-1}\quad\alpha\in\mathbb{D}_{\rho}~.
\end{equation}
With the same $\alpha$,
let $\psi_{\alpha}(z)$ be the scalar unitary function,
\begin{equation}\label{eq:ScalarRational}
\psi_{\alpha}(z^N):=\frac{1-|{\alpha}|^2}{z^N-\alpha}=
c(zI_N-A)^{-1}b,
\end{equation}
for some $A\in\C^{N\times N}$, $b\in\C^{N\times 1}$
and $c\in\C^{1\times N}$.

Then the $2N\times 2N$ system matrix $M$ associated
with $V(z^N)$ in \eqref{eq:u(z^N)} is given by,
\[
M=\left(\begin{array}{r|r}
A&bv^*\\ \hline vc &I_N-(1+{\alpha}^*)vv^*
\end{array}\right).
\]
\end{Ob}
\vskip 0.2cm

Indeed, it is straightforward to verify that one can equivalently
write $V(z^N)$ in \eqref{eq:u(z^N)} as,
\[
V(z^N)=v\psi_{\alpha}(z)v^*+I_N-(1+{\alpha}^*)vv^*=
vc(zI_N-A)^{-1}bv^*+I_N-(1+{\alpha}^*)vv^*.
\]
This is illustrated next.

\begin{Ex}\label{Ex:U(z^2)}
{\rm

1. Consider the system matrix $\tilde{M}_{\alpha}$
associated with the scalar function $\psi_{\alpha}(z^N)$
in \eqref{eq:ScalarRational} with $N=2$. It depends on the
parameter $\alpha\in\mathbb{D}_{\rho}$
and is given by,
\[
\tilde{M}_{\alpha,2}=\left(\begin
{array}{cr|c}
\sqrt{\alpha}~~&1
&0\\ 0~~&-\sqrt{\alpha}
&\sqrt{1-|{\alpha}|^2}\\
\hline
\sqrt{1-|{\alpha}|^2}~~&~0
&0\end{array}\right)~.
\]
\vskip 0.2cm

2. Now, using Observation \ref{Ob:MtildeU} with $v=${\mbox{\tiny$
\begin{pmatrix}0\\1\end{pmatrix}$}}, the system matrix
associated with the $2\times 2$ valued $V_a(z^2)$
in \eqref{eq:Ua} is again depending on $\alpha$ and given by,
\begin{equation}\label{eq:RealizationU(z^2alpha)}
M_{\alpha}=\left(\begin{array}{cr|rc}
\sqrt{\alpha}~~&1&0&0\\
0~~&-\sqrt{\alpha}&0&\sqrt{1-|{\alpha}|^2}\\
\hline
0~~&~~0&1&0\\ \sqrt{1-|{\alpha}|^2}~~&~0
&0&-{\alpha}^*\end{array}\right)~.
\end{equation}
\vskip 0.2cm

3. Now consider the system matrix $\tilde{M}_{\alpha}$
associated with the scalar function $\psi_{\alpha}(z^N)$
in \eqref{eq:ScalarRational} with $N=4$. It depends on the
parameter $\alpha\in\mathbb{D}_{\rho}$ and is given by,
\[
\tilde{M}_{\alpha,4}=\left(
\begin{array}{crrr|c}
\alpha^{\frac{1}{4}}&1&~~~0&0&0\\
0~~&-\alpha^{\frac{1}{4}}&~1&~0~&0\\
0~~&0&i\alpha^{\frac{1}{4}}&1&0\\
0~~&0&~~~0&-i\alpha^{\frac{1}{4}}&\sqrt{1-|\alpha |^2}\\
\hline
\sqrt{1-|\alpha |^2}&0&0&0&0\end{array}\right).
\]
4. Next, substituting $\alpha=\beta$ in the above $\tilde{M}_{\alpha,4}$
and using Observation \ref{Ob:MtildeU} with
$v=${\mbox{\tiny$\begin{pmatrix}1\\0\end{pmatrix}$}}, yields the system
matrix associated with $V_{\beta}(z^2)$ in \eqref{eq:Ub}. It
depends on the parameter $\beta\in\mathbb{D}_{\rho}$ and is given by,
\begin{equation}\label{eq:RealizationU(z^4beta)}
M_{\beta,4}=\left(\begin{array}{crrr|cr}
\beta^{\frac{1}{4}}&1&~~~0&0&0&0\\
0~~&-\beta^{\frac{1}{4}}&~1&~0~&0&0\\
0~~&0&i\beta^{\frac{1}{4}}&1&0&0\\
0~~&0&~~~0&-i\beta^{\frac{1}{4}}
&\sqrt{1-|\beta |^2}&0\\
\hline
\sqrt{1-|\beta |^2}&0&0&0&-\beta^*&0\\
0&0&0&0&~0&~1\end{array}\right).
\end{equation}
}\qed
\end{Ex}

\section{A recipe for state space realization of
wavelet filters}
\label{sec:8} \setcounter{equation}{0}

Recall that in Proposition \ref{Pr:WF} we introduced the following
construction scheme of wavelet filters. By multiplying from
the left a given \mbox{$W_a\in\mathcal{W}_N$}, of McMillan degree
$N\left(\frac{1}{2}(N-1)+m\right)$, by an elementary unitary matrix
$V(z^N)$, one can construct another wavelet filter
$W_b(z)$ of McMillan degree \mbox{$N\left(\frac{1}{2}(N+1)+m\right)$.}
Indeed,
\[
W_b(z)=V(z^N)W_a(z),
\]
where $V\in\mathcal{U}$ is elementary with parameters
\mbox{$v\in\mathbb{S}_{N-1}$} and
\mbox{$\alpha\in\mathbb{D}_{\rho}$,} see \eqref{eq:tildeU}.
\vskip 0.2cm

Assuming the state space realization of
$W_a(z)$ is already known, the state space realization of
$V(z^N)$ is given in Example \ref{Ex:U(z^2)}
and thus $M_b$, the system matrix associated with $W_b(z)$
can be obtained through Observation \ref{Ob:series}.

A somewhat similar idea appeared in \cite[Section 4]{Wa}.

\begin{Ex}\label{Ex:RealizationWaWb}
{\rm
We here construct, in two steps, the system matrix $M_b$
associated with $W_b(z)$ in \eqref{eq:deg7filter}, a
wavelet filter of dimension $N=2$ and index $m=3$.
\vskip 0.2cm

Recall that $W_a(z)$ in \eqref{eq:deg3filter} is of
McMillan degree three. Using Observation \ref{Ob:series},
from the realizations in item 1 of Example \ref{Ex:canonicalM}
and in \eqref{eq:RealizationU(z^2alpha)} (item 2 in Example
\ref{Ex:U(z^2)}), we construct $M_a$, the system matrix
associated with $W_a(z)$,
\begin{equation}\label{eq:RealizationWa}
M_a=\left(\begin
{array}{rcr|rr}
\sqrt{\alpha}~~&1-|\alpha|^2&0&0&~~0\\
0~~&-\sqrt{\alpha}&1&0&~~0\\
0~~&0&0&\frac{1}{\sqrt{2}}&-\frac{1}{\sqrt{2}}\\
\hline
0~~&0&0
&\frac{1}{\sqrt{2}}&~\frac{1}{\sqrt{2}}\\
1~~&0&-\alpha^*&0&~~0\end{array}\right)~.
\end{equation}
\vskip 0.2cm

Next, recall that $W_b(z)$ in \eqref{eq:deg7filter} is of
McMillan degree seven. Using Observation \ref{Ob:series},
again, from the realizations in \eqref{eq:RealizationWa} and in
\eqref{eq:RealizationU(z^4beta)} (item 3 in Example \ref{Ex:U(z^2)}),
we construct $M_b$, the system matrix associated with
$W_b(z)$,
\begin{equation}\label{eq:Mb}
M_b=\left(\begin{array}{crrrcrr|cc}
\beta^{\frac{1}{4}}&1&~~~0&0&0&0~~&0&0&0\\
0~~&-\beta^{\frac{1}{4}}&~1&~0~&0&0~~&0&0&0\\
0~~&0&i\beta^{\frac{1}{4}}&1&0&0~~&0&0&0\\
0~~&0&~~~0&-i\beta^{\frac{1}{4}}&0&0~~&0
&\frac{\sqrt{1-|\beta |^2}}{\sqrt{2}}&
\frac{\sqrt{1-|\beta |^2}}{\sqrt{2}}\\
0~~&0&~~~0&0&\sqrt{\alpha}~~&1&~~0&~0&~0\\
0~~&0&~~~0&0&0&-\sqrt{\alpha}~~&1-|\alpha|^2&~0&~0\\
0~~&0&~~~0&~~0&~0&0&0
&~\frac{1}{\sqrt{2}}&-\frac{1}{\sqrt{2}}\\
\hline
\sqrt{1-|\beta |^2}&0&0&0&~0~&0~~&0
&-\frac{\beta^*}{\sqrt{2}}&-\frac{\beta^*}{\sqrt{2}}\\
0&0&0&0&\sqrt{1-|\alpha |^2}~&0~~&-\alpha^*
&~0&~0\end{array}\right).
\end{equation}
}\qed \end{Ex}

\section{Concluding remarks}
\label{sec:10}
\setcounter{equation}{0}

1. The parameterization of all wavelet filters in Observation
\ref{Ob:parameters} has several advantages:
\begin{itemize}

\item One can use the index $m$ as a design parameter.

\item One has a clear view on the trade-off involved in
increasing the index $m$: Additional degrees of freedom on
the expense of computational burden, which is even more
apparent from the state-space model.

\item The spectral radius $\rho$, is another design
parameter. This point is extended in item 4 below.

\item $\rho$ offers another the trade-off.
For $1>>\rho>0$ one obtains ``almost FIR" filters,
which is a much larger family than FIR's. However,
for prescribed $m$, increasing $\rho$, weakens
the noise attenuation and the slows down the
rate of convergence.

\item The set of design parameters is convex and thus
convenient for design through optimization, see e.g.
\cite{CZZM4}, \cite{HHN}, \cite{VHEK} and \cite{TV}.
\end{itemize}
\vskip 0.2cm

2. The stage-by-stage construction of state space realization
introduced in Sections \ref{sec:7}, \ref{sec:8} is not only
easy to employ, but also reveals the appealing structure of the
realization of wavelet filters, see e.g. \eqref{eq:Mb}. This is
significant in facilitating the computations involved when one
wishes to employ Matlab based algorithms such as LMI toolbox
as proposed in e.g. \cite{CZZM2}, \cite{CZZM4}, \cite{VHEK},
or solvers of the Riccati equation, see e.g.
\cite{CZZM1}, \cite{CZZM3}.
\vskip 0.2cm

3. Further properties of rational wavelet filters

In the present paper we concentrated on a systems-theoretic
approach to a parameterization of the family of all wavelet
filters as matrix functions. There are important applications
which we only sketch briefly here, for example:\\
(a) Within the entire family, identify and analyze those wavelet
filters that have vanishing moments; (the vanishing-moment
wavelets have especially attractive algorithms for approximation!);\\
(b) Give numerical recipes for how to translate our parameterization
into wavelet functions on the real line or on more general ambient
spaces, so obtaining scaling functions, wavelet generators
(father/mother functions) and associated wavelet bases from our
filter matrix-functions.\\
Both (a) and (b) are major topics in wavelet theory, see e.g.
\cite{BFMP09, BlKr12, GlDA09,  Gu11,   Ha11,
 LaRa06,  MiSu07, OlOl10,
Rae09,SuVi10, YaZh08}.
\\

4. Relaxing Schur stability (taking $\rho>1$)

In the description of the set $\mathcal{U}$ in \eqref{eq:unit}
we confined the discussion to Schur asymptotically stable
functions, i.e. analytic outside $\mathbb{D}_{\rho}$
for some prescribed $\rho\in[0,~1]$.
Consequently, in the stage-by-stage construction of
$\mathcal{U}$ in Section \ref{sec:3}, we take the poles to be
within $\mathbb{D}_{\rho}$~.
\vskip 0.2cm

However, in some applications this requirement is not imperative,
e.g. where the filters are not
associated with dynamical systems. Then, one can extend the
discussion and allow the poles of a function in $\mathcal{U}$
to be in $\mathbb{D}_{\rho}\smallsetminus{\mathbb T}$, for
some $\rho>1$, see e.g. \cite[Section 3]{AlGo1} and in the
framework of wavelet filters see \cite{AJL}, \cite{AJLM}.
\vskip 0.2cm

5. Not necessarily rational wavelet filters.

Equation \eqref{eq:tildeU} may be substituted by,
\[
\tilde{V}(z^N):=I_N+\left(\psi(z^N)
-1\right)vv^*,
\]
where $~v\in\mathbb{S}_{N-1}$ is a parameter and the scalar
function $\psi(z)$, mapping $\mathbb{T}$ to itself, is
with the appropriate analyticity conditions, e.g.
meromorphic outside $\mathbb{D}$. For example, one
can take
\[
\psi(z^N)=e^{\frac{z^N+1}{z^N-1}},
\]
which maps the open unit disk to itself.

Thus, if $W(z)$ is a wavelet filter satisfying
\eqref{eq:WaveletFilters}, then so is the product
\[
\tilde{V}(z^N)W(z).
\]
A thorough study of non rational generalizations of wavelet
filters can be found in \cite{AJL}.
\vskip 0.2cm

6. Non-square wavelet filters

The property $\mathcal{C}$ defined in \eqref{eq:m1} can in fact
be generalized to having
\[
F(z)=\begin{pmatrix} \hat{f}_0(z)&
\hat{f}_0(\e z)&\cdots &\hat{f}_0(\e^{N-1}z)\\
\hat{f}_1(z)&\hat{f}_1(\e z)&\cdots &\hat{f}_1(\e^{N-1}z)\\
\vdots& &  & \\
\hat{f}_{l-1}(z)&\hat{f}_{l-1}(\e z)&\cdots &
\hat{f}_{l-1}(\e^{N-1}z)
\end{pmatrix},
\]
satisfying the condition \eqref{eq:sym1},
\[
F(\e z)=F(z)\hat{P},
\]
where $\hat{P}$ is as before.
\vskip 0.2cm

Now for the requirement on $\mathbb{T}$ we return to \eqref{RE}.
For $N\geq l$, this requirement is generalized
to coisometry $FF^*=I_l$.

If however $l>N$ (in signal processing literature this is
referred to as the ``oversampling" case. In other places
``frames" are discussed) we return to \eqref{RE} and take,
\[
An(z)=F(z)\quad\quad\quad
Sy(z)=\left(F\left(\frac{1}{z^*}\right)^*F(z)\right)^{-1}
\cdot{F\left(\frac{1}{z^*}\right)^*}.
\]
In particular, $Sy(z)$ is a left inverse of $An(z)$, on $\mathbb{T}$.
For mathematical analysis of this case see e.g. \cite{AJL}, \cite{J5}.
From engineering point of view this case is addressed in a series of
papers \cite{CZZM1}, \cite{CZZM2}, \cite{CZZM3}, \cite{CZZM4}.
\vskip 0.2cm

7. Relaxing the condition \eqref{RE} on all $\mathbb{T}$.

One can require that $SyAn\approx I$ only on part of
$\mathbb{T}$. This leads to ``wavelets on fractals", see
e.g.  \cite{DJ1}, \cite{DJ2}, \cite{DJ3}. See also \cite{HHN}.
\vskip 0.2cm

\section{Appendix -additional background}
\label{sec:11}
\setcounter{equation}{0}

\subsection{The $Z$-transform of $N$-decimated $N$-expanded
sequence}\label{Ap1}
\vskip 0.2cm

We next address the passage from $G(z)$ to $G(z^N)$ and recall
in the following background. Let $\{ a_k\}_{k=-\infty}^{\infty}$
be a sequence. For a natural $N$, its {\em decimated} version is
the sequence, $\{ b_k\}_{k=-\infty}^{\infty}$ with $b_k=a_{Nk}$.
The $N$-{\em expanded} version of $\{ b_k\}_{k=-\infty}^{\infty}$
is the sequence  $\{ c_k\}_{k=-\infty}^{\infty}$ where
all elements are zero, except $c_{Nk}=b_k$. Alternatively,
one can write that all elements in $\{ c_k\}_{k=-\infty}^{\infty}$
are zero, except $c_{Nk}=a_{Nk}$, see e.g. \cite[Chapter 3]{SN},
\cite[Subsection 4.1.1]{Va}.
\vskip 0.2cm

Now, we recall that if the $Z$-transform of
$\{ a_k\}_{k=-\infty}^{\infty}$ exists and denoted by $F(z)$,
the $Z$-transform of $\{ c_k\}_{k=-\infty}^{\infty}$ exists
and is given by $F(z^N)$.
Indeed, assume that $\{ a_k\}_{k=-\infty}^{\infty}$ is
so that all elements to the left of $a_{-k_o}$, for some
positive finite $k_o$, vanish. Thus, $\{ a_k\}$ can be
written as a sum $A_-+A_+$ with $A_-=\{ a_k\}_{k=-k_o}^0$
and $A_+=\{ a_k\}_1^{\infty}$. Clearly, the $Z$-transform of
$A_-$, denoted by $F_-(z)$ is well defined and is rational.
Assume in addition that the $Z$-transform of $A_+$,
denoted by $F_+(z)$, is rational. Then, the
fact that the $Z$-transform of $\{c_k\}_{k=-k_o}^{0}$ exists
and is rational, is obvious. Now to see that also
$Z$-transform of $\{c_k\}_{k=1}^{\infty}$ exists and is
rational we examine the associated Hankel matrices.
Let $H_a$ and $H_c$ be the Hankel matrices associated
with $\{ a_k\}_{k=1}^{\infty}$ and
$\{c_k\}_{k=1}^{\infty}$, respectively, e.g.
\[
H_a=\left(\begin{smallmatrix}
a_{11}&a_{12}&a_{13}&a_{14}&a_{15}&a_{16}&a_{17}\\
a_{12}&a_{13}&a_{14}&a_{15}&a_{16}&a_{17}&a_{28}\\
a_{13}&a_{14}&a_{15}&a_{16}&a_{17}&a_{28}&a_{38}\\
a_{14}&a_{15}&a_{16}&a_{17}&a_{28}&a_{38}&a_{48}\\
a_{15}&a_{16}&a_{17}&a_{28}&a_{38}&a_{48}&a_{58}\\
a_{16}&a_{17}&a_{28}&a_{38}&a_{48}&a_{58}&a_{68}\\
a_{17}&a_{28}&a_{38}&a_{48}&a_{58}&a_{68}&a_{78}
\end{smallmatrix}\right)\quad\quad
H_c=\left(\begin{smallmatrix}
a_{11}&0&0&a_{14}&0&0&a_{17}\\
0&0&a_{14}&0&0&a_{17}&0\\
0&a_{14}&0&0&a_{17}&0&0\\
a_{14}&0&0&a_{17}&0&0&a_{48}\\
0&0&a_{17}&0&0&a_{48}&0\\
0&a_{17}&0&0&a_{48}&0&0\\
a_{17}&0&0&a_{48}&0&0&a_{78}
\end{smallmatrix}\right).
\]
Taking
\[
E:=\left(\begin{smallmatrix}
1&0&0&0&0&0&0\\
0&0&0&1&0&0&0\\
0&0&0&0&0&0&1\end{smallmatrix}\right),
\]
one can always write,
\[
\hat{h}=EH_aE^*=EH_cE^*\quad\quad\quad
H_c=E^*\hat{h}E\quad\quad\quad
\hat{h}=\left(\begin{smallmatrix}
a_{11}&a_{14}&a_{17}\\
a_{14}&a_{17}&a_{48}\\
a_{17}&a_{48}&a_{78}\end{smallmatrix}\right).
\]
Namely, $\hat{h}$ is another Hankel matrix. This implies
that
\[
{\rm rank}(H_c)={\rm rank}(E^*\hat{h}E)=
{\rm rank}(\hat{h})\leq {\rm rank}(H_a).
\]
By assumption the rank of $H_a$ is finite then so is
${\rm rank}(H_c)$, which implies that both
$\{ a_k\}_{k=1}^{\infty}$ and $\{c_k\}_{k=1}^{\infty}$ admit
state space realization, see e.g. \cite[Lemma 6.5-7]{Ka},
\cite[Lemma 5.5.5]{So} and thus the $Z$-transform of each
sequence exists and is rational.
\vskip 0.2cm

\subsection{Functions in $\mathcal{U}$ admitting
state space realization}\label{Ap3}

We first cites
an adapted version of \cite[Theorem 3.9]{AlGo1} (a closely
related result appeared in \cite[Theorem 3]{GVKDM})

\begin{Tm}\label{Tm:SteinRealizationU}
Let $U(z)$ analytic at infinity, be of McMillan degree $p$, and let
\[
U(z)=C(zI-A)^{-1}B+D
\]
be a minimal realization of $U$. Then,
\[
U\in\mathcal{U}_N
\]
if and only if there exists a $p\times p$ non-singular
Hermitian matrix $H$ (uniquely
determined from the given realization) such that
\[
\left(\begin{smallmatrix}
A&B\\ C&D\end{smallmatrix}\right)^*
\left(\begin{smallmatrix}
H&0\\0&I_N\end{smallmatrix}\right)
\left(\begin{smallmatrix} A&B\\C&D\end{smallmatrix}\right)=
\left(\begin{smallmatrix}
H&0\\0&I_N\end{smallmatrix}\right).
\]
\end{Tm}

This $H$ is called the {\sl Hermitian matrix associated}
(with the given minimal realization).
\vskip 0.2cm

Using the system matrix $M$ notation from \eqref{eq:SystemMatrix},
the Stein Equation in Theorem \ref {Tm:SteinRealizationU} can be
compactly written as
\[
M^*\hat{H}M=\hat{H}\quad\quad\quad
\hat{H}:={\rm diag}\{H,~I_N\}.
\]
\begin{center}
Acknowledgement
\end{center}

The authors wish to thank the referees for helping improving
this work.

\bibliographystyle{plain}

\end{document}